\documentclass[12pt, a4paper]{amsart}
\usepackage[utf8]{inputenc}
\usepackage{amsxtra,amssymb,amsthm,amsmath,amscd,mathrsfs, epsfig, eufrak}
\usepackage{amscd, amsmath, mathrsfs, amssymb, amsthm, amsxtra, bbding, epsfig, eucal, eufrak, graphicx, latexsym, mathrsfs, mathbbol, bbold}
\usepackage[all]{xy}

\usepackage{tikz}

\usepackage[normalem]{ulem}
\usepackage{soul}
\usepackage{color}

\setstcolor{red}

\makeatletter
\@namedef{subjclassname@2010}{%
\textup{2010} Mathematics Subject Classification}
\makeatother

\def \N {{\mathbb N}}
\def \P {{\mathbb P}}

\def \d {\,{\rm d}}

\def\Li{\hbox{{\rm Li}}}

\def\le{\leqslant}
\def\ge{\geqslant}

\theoremstyle{plain}
\newtheorem{theorem}{Theorem}

\theoremstyle{remark}

\numberwithin{equation}{section}

\addtocounter{footnote}{1}

\begin{document}


\title[\tiny On the primes in floor function sets]
{On the primes in floor function sets}
\author[\tiny Rong Ma \& Jie Wu]{Rong Ma \& Jie Wu}

\address{%
Rong Ma
\\
School of Mathematics and Statistics
\\
Northwestern Polytechnical University
\\
Xi'an
\\
Shaanxi 710072
\\
China}
\email{marong0109@163.com}

\address{%
Jie Wu
\\
CNRS UMR 8050
\\
Laboratoire d'Analyse et de Math\'ematiques Appliqu\'ees
\\
Universit\'e Paris-Est Cr\'eteil
\\
94010 Cr\'eteil cedex
\\
France}
\email{jie.wu@u-pec.fr}

\date{\today}

\begin{abstract}
Let $[t]$ be the integral part of the real number $t$ and 
let $\mathbb{1}_{\P}$ be the characteristic function of the primes.
Denote by $\pi_{\mathcal{S}}(x)$ the number of primes 
in the floor function set $\mathcal{S}(x) := \{[\frac{x}{n}] : 1\le n\le x\}$
and by $S_{\mathbb{1}_{\P}}(x)$ the number of primes in the sequence $\{[\frac{x}{n}]\}_{n\ge 1}$.
Very recently, Heyman proves
$$
\pi_{\mathcal{S}}(x) = \frac{4\sqrt{x}}{\log x} + O\bigg(\frac{\sqrt{x}}{(\log x)^2}\bigg),
\qquad
S_{\mathbb{1}_{\P}}(x) = C_{\mathbb{1}_{\P}} x + O(x^{1/2})
$$
for $x\to\infty$, where $C_{\mathbb{1}_{\P}} := \sum_{p} \frac{1}{p(p+1)}$.
In this short note, we propose better results
$$
\pi_{\mathcal{S}}(x) 
= \int_2^{\sqrt{x}} \frac{\d t}{\log t} + \int_2^{\sqrt{x}} \frac{\d t}{\log(x/t)}
+ O\Big(\sqrt{x}\,{\rm e}^{-c(\log x)^{3/5}(\log\log x)^{-1/5}}\Big),
$$
and
$$
S_{\mathbb{1}_{\P}}(x)
= C_{\mathbb{1}_{\P}} x + O_{\varepsilon}(x^{9/19+\varepsilon})
$$
for $x\to\infty$, 
where $c>0$ is a positive constant and $\varepsilon$ is an arbitrarily small positive number.
\end{abstract}


\maketitle

\section{Introduction}

As usual, denote by $\pi(x)$ the number of primes $p\le x$.
It is well-known that
\begin{itemize}
\item[(a)]
The prime number theorem states as follows
\begin{equation}\label{PNT:weak}
\pi(x) = \frac{x}{\log x} + O\bigg(\frac{x}{(\log x)^2}\bigg)
\qquad
(x\to\infty).
\end{equation}

\item[(b)]
A strong form of this theorem is the following
\begin{equation}\label{PNT:strong}
\pi(x) = \Li(x)  + O(x \exp(-c(\log x)^{3/5}(\log_2x)^{-1/5}))
\end{equation}
for $x\to\infty$,
where $c$ is a positive constant, $\log_2$ denotes the iterated logarithm function and
$$
\Li(x) := \int_2^x \frac{\d t}{\log t}\cdot
$$

\item[(c)]
The Riemann hypothesis is equivalent to the asymptotic formula
\begin{equation}\label{PNT:RH}
\pi(x) = \Li(x)  + O_{\varepsilon}(x^{1/2+\varepsilon})
\qquad\;
(x\to\infty),
\end{equation}
where $\varepsilon$ is an arbitrarily small positive number.
\end{itemize}

Let $[t]$ be the integral part of the real number $t$.
Recently, Bordell\`es, Dai, Heyman, Pan and Shparlinski \cite{BDHPS2019} 
proposed to investigate the asymptotic behaviour of summative function
\begin{equation}\label{def:Sfx}
S_f(x)
:= \sum_{n\le x} f\Big(\Big[\frac{x}{n}\Big]\Big)
\end{equation}
under some simple hypothesis on the growth of $f$,
and this problem has received attention of many authors 
\cite{Wu2020, Zhai2020, Bordelles2020, LiuWuYang2021a, LiuWuYang2021b}.
If we use $\Lambda(n)$ to denote the von Mangoldt function, 
then \cite[Theorem 1.2(i)]{Wu2020} or \cite[Theorem 1]{Zhai2020} give us immediately
\begin{equation}\label{Wu-Zhai:1}
S_{\Lambda}(x) = C_{\Lambda} x + O_{\varepsilon}(x^{1/2+\varepsilon}),
\end{equation}
for any $\varepsilon>0$ and $x\to\infty$, where $C_{\Lambda} := \sum_{n\ge 1} \frac{\Lambda(n)}{n(n+1)}$.
However Ma and Wu \cite{MaWu2020} applied the Vaughan identity 
and the technique of one-dimensional exponential sums to break the $\frac{1}{2}$-barrier by establishing
$$
S_{\Lambda}(x) = C_{\Lambda} x + O_{\varepsilon}(x^{35/71+\varepsilon}).
$$
This result seems rather interesting if we compare it with the assertion (c) above.
The exponent $\frac{35}{71}$ has been improved to $\frac{97}{203}$ by Bordell\`es \cite{Bordelles2020} 
and $\frac{9}{19}$ by Liu-Wu-Yang \cite{LiuWuYang2021a}, respectively,
with the help of more sophistic technique of exponential sums.

Let $\P$ be the set of all primes and let ${\P}_{\rm ower}$ be the set of all prime powers.
Denote by $\mathbb{1}_{\P}$ and $\mathbb{1}_{{\P}_{\rm ower}}$ their characteristic functions, respectively.
Very recently Heyman \cite{Heyman2021} proposed to study the number of primes 
in the floor function set $\mathcal{S}(x) := \{[\frac{x}{n}] : 1\le n\le x\}$:
\begin{equation}\label{def:piGx}
\pi_{\mathcal{S}}(x)
:= \sum_{\substack{p\le x\\ \exists\,n\in \N\,\text{such that}\; [x/n]=p}} 1
\end{equation}
and the number of primes or prime powers in the sequence $\{[\frac{x}{n}]\}_{n\ge 1}$:
\begin{equation}\label{def:S1Px}
S_{\mathbb{1}_{\P}}(x)
:= \sum_{n\le x} \mathbb{1}_{\P}\Big(\Big[\frac{x}{n}\Big]\Big),
\qquad
S_{\mathbb{1}_{{\P}_{\rm ower}}}(x)
:= \sum_{n\le x} \mathbb{1}_{{\P}_{\rm ower}}\Big(\Big[\frac{x}{n}\Big]\Big).
\end{equation}
The principal result of Heyman \cite[Theorem 1]{Heyman2021} is the following asymptotic formula 
\begin{equation}\label{Heyman:1}
\pi_{\mathcal{S}}(x) = \frac{4\sqrt{x}}{\log x} + O\bigg(\frac{\sqrt{x}}{(\log x)^2}\bigg)
\end{equation}
as $x\to\infty$.
This is the prime number theorem in weak form for the set $\mathcal{S}(x)$,
i.e. analogue of \eqref{PNT:weak} for this set.
This result is rather interesting, since $\mathcal{S}(x)$ is a very spare subset of $[1, x]\cap \N$.
In fact Heyman \cite[Theorems 1 and 2]{Heyman2019} has been proved that
\begin{equation}\label{Heyman:result}
|\mathcal{S}(x)| = 2\sqrt{x} + O(1)
\end{equation}
for $x\to\infty$.
Probably this is the first example of a such spare subset of $[1, x]\cap \N$ 
for which the prime number theorem holds.

It seems natural and interesting to establish analogue of \eqref{PNT:strong}, i.e.
the prime number theorem in strong form for the set $\mathcal{S}(x)$.
The first aim of this short note is to prove a such result.

\begin{theorem}\label{thm1}
{\rm (i)}
For $x\to\infty$, we have
\begin{equation}\label{thm1:eq1}
\pi_{\mathcal{S}}(x)
= \Li_{\mathcal{S}}(x) + O(\sqrt{x} \exp(-c'(\log x)^{3/5}(\log_2x)^{-1/5})),
\end{equation}
where $c'>0$ is a positive constant and
\begin{equation}\label{thm1:eq2}
\Li_{\mathcal{S}}(x)
:= \int_2^{\sqrt{x}} \frac{\d t}{\log t} + \int_2^{\sqrt{x}} \frac{\d t}{\log(x/t)}\cdot
\end{equation}

{\rm (ii)}
There is a real sequence $\{a_n\}_{n\ge 1}$ with $a_1=4$ such that for any positive integer $N\ge 1$ we have
\begin{equation}\label{thm1:eq3}
\pi_{\mathcal{S}}(x)
= \sqrt{x} \sum_{n=1}^{N} \frac{a_n}{(\log x)^n} + O_N\bigg(\frac{\sqrt{x}}{(\log x)^{N+1}}\bigg)
\end{equation}
as $x\to\infty$.
\end{theorem}

Theorems 5 and 7 of \cite{Heyman2021} can be stated as follows:
\begin{align}
S_{\mathbb{1}_{\P}}(x)
& = C_{\mathbb{1}_{\P}} x + O(x^{1/2}),
\label{Heyman:2}
\\
S_{\mathbb{1}_{{\P}_{\rm ower}}}(x)
& = C_{\mathbb{1}_{{\P}_{\rm ower}}} x + O(x^{1/2}),
\label{Heyman:3}
\end{align}
where $C_{\mathbb{1}_{\P}} := \sum_{p} \frac{1}{p(p+1)}$ and
$C_{\mathbb{1}_{{\P}_{\rm ower}}} := \sum_{p, \, \nu\ge 1} \frac{1}{p^{\nu}(p^{\nu}+1)}$.
Similar to \eqref{Wu-Zhai:1}, 
these are immediate consequences of \cite[Theorem 1.2(i)]{Wu2020} or \cite[Theorem 1]{Zhai2020}.
Heyman \cite[Theorem 6]{Heyman2021} also proved that there is a positive constant $A>0$ such that the inequality
\begin{equation}\label{Heyman:4}
S_{\mathbb{1}_{\P}}(x)
\ge C_{\mathbb{1}_{\P}} x - Ax^{1/2}/\log x
\end{equation}
holds for $x\ge 2$.

The second aim of this note is to propose better results by breaking the $\frac{1}{2}$-barrier in the error term.

\begin{theorem}\label{thm2}
For any $\varepsilon>0$, we have
\begin{align}
S_{\mathbb{1}_{\P}}(x)
& = C_{\mathbb{1}_{\P}} x + O_{\varepsilon}(x^{9/19+\varepsilon}),
\label{thm2:eq1}
\\
S_{\mathbb{1}_{{\P}_{\rm ower}}}(x)
& = C_{\mathbb{1}_{{\P}_{\rm ower}}} x + O_{\varepsilon}(x^{9/19+\varepsilon}),
\label{thm2:eq2}
\end{align}
as $x\to\infty$, where the implied constants depend on $\varepsilon$.
\end{theorem}

We note that
very recently Yu and Wu \cite{YuWu2022} generalised Heyman's \eqref{Heyman:result} by showing
\begin{equation}\label{YuWu:result}
\mathcal{S}(x; q, a)
:= \sum_{\substack{m\in \mathcal{S}(x)\\ m\equiv a ({\rm mod}\,q)}} 1
= \frac{2\sqrt{x}}{q} + O((x/q)^{1/3}\log x)
\end{equation}
uniformly for $x\ge 3$, $1\le q\le x^{1/4}/(\log x)^{3/2}$ and $1\le a\le q$,
where the implied constant is absolute.
This confirms a recent numeric test of Heyman.

\vskip 8mm

\section{Proof of Theorem \ref{thm1}}

We begin by following the argument of \cite{Heyman2021}.
Firstly we note that
$$
\mathcal{S}(x) 
= \Big\{p\in \mathbb{P} : \exists \; n\in [1, x] \;\text{such that}\; \Big[\frac{x}{n}\Big]=p\Big\}.
$$ 
Further, if $\big[\frac{x}{n}\big] = p\in \mathbb{P}$, then $x/(p+1)<n\le x/p$. 
Thus we can write
\begin{equation}\label{G:expression}
\pi_{\mathcal{S}}(x) = \sum_{p\le x} \mathbb{1}\Big(\Big[\frac{x}{p}\Big]-\Big[\frac{x}{p+1}\Big]>0\Big)
= G_1(x) + G_2(x),
\end{equation} 
where $\mathbb{1}=1$ if the statement is true and 0 otherwise, and
\begin{align*}
G_1(x)
& := \sum_{p\le \sqrt{x}} \mathbb{1}\Big(\Big[\frac{x}{p}\Big]-\Big[\frac{x}{p+1}\Big]>0\Big),
\\
G_2(x)
& := \sum_{\sqrt{x}<p\le x} \mathbb{1}\Big(\Big[\frac{x}{p}\Big]-\Big[\frac{x}{p+1}\Big]>0\Big).
\end{align*}

For $p\le \sqrt{x}-1$, we have
$$
\Big[\frac{x}{p}\Big]-\Big[\frac{x}{p+1}\Big]
>\frac{x}{p(p+1)}-1
>0.
$$
Thus the prime number theorem \eqref{PNT:strong} gives us
\begin{equation}\label{G1}
\begin{aligned}
G_1(x) 
& = \pi(\sqrt{x}) + O(1)
\\
& = \Li(\sqrt{x})  + O(x \exp(-c'(\log x)^{3/5}(\log_2x)^{-1/5}))
\end{aligned}
\end{equation}
for $x\ge 3$, where $c'>0$ is a positive constant.

Next we treat $G_2(x)$.
Noticing that
$$
0<\frac{x}{p}-\frac{x}{p+1}=\frac{x}{p(p+1)}<1
$$
for $p>\sqrt{x}$, the quantity $\big[\frac{x}{p}\big]-\big[\frac{x}{p+1}\big]$ can only equal to 0 or 1. 
On the other hand, for $p>x^{10/19}$, then $p=[\frac{x}{n}]$ for some $n\le x^{9/19}$.
Thus we can write
\begin{equation}\label{G2:decomposition}
\begin{aligned}
G_2(x)
& = \sum_{x^{1/2}<p\le x^{10/19}} \Big(\Big[\frac{x}{p}\Big]-\Big[\frac{x}{p+1}\Big]\Big) + O(x^{9/19})
\\
& = \sum_{x^{1/2}<p\le x^{10/19}} \Big(\frac{x}{p}-\frac{x}{p+1}-\psi\Big(\frac{x}{p}\Big)+\psi\Big(\frac{x}{p+1}\Big)\Big) + O(x^{9/19})
\\
& = G_{2, 1}(x) - G_{2, 2}^{\langle 0\rangle}(x) + G_{2, 2}^{\langle 1\rangle}(x) + O(x^{9/19}),
\end{aligned}
\end{equation}
where $\psi(t):=t-[t]-\frac{1}{2}$ and
\begin{align*}
G_{2, 1}(x)
& := \sum_{x^{1/2}<p\le x^{10/19}} \Big(\frac{x}{p}-\frac{x}{p+1}\Big),
\\
G_{2, 2}^{\langle\delta\rangle}(x)
& := \sum_{x^{1/2}<p\le x^{10/19}} \psi\Big(\frac{x}{p+\delta}\Big)
\quad
(\delta=0, 1).
\end{align*}
With the help of the prime number theorem \eqref{PNT:strong}, a simple partial integration allows us to derive that
\begin{align*}
G_{2, 1}(x)
& = \sum_{x^{1/2}<p\le x/2} \frac{x}{p^2} + O(x^{9/19})
= x \int_{\sqrt{x}}^{x/2} \frac{\d \pi(t)}{t^2} + O(x^{9/19})
\\
& = x \int_{\sqrt{x}}^{x/2} \frac{\d t}{t^2\log t}
+ O\big(\sqrt{x}\exp(-c'(\log x)^{3/5}(\log_2x)^{-1/5}\big),
\end{align*}
where $c'>0$ is a positive constant.
Making the changement of variables $t\to x/t$ in the last integral, it follows that
\begin{equation}\label{G21}
G_{2, 1}(x)
= \int_2^{\sqrt{x}} \frac{\d t}{\log(x/t)}
+ O\big(\sqrt{x}\exp(-c'(\log x)^{3/5}(\log_2x)^{-1/5}\big)
\end{equation}
for $x\to\infty$.

It remains to bound $G_{2, 2}^{\langle\delta\rangle}(x)$.
Similar to \cite{LiuWuYang2021a}, define
$$
\mathfrak{S}_{\delta}(x; D, D')
:= \sum_{D<d\le D'} \Lambda(d) \psi\Big(\frac{x}{d+\delta}\Big).
$$
According to \cite[(4.3)]{LiuWuYang2021a}, for any $\varepsilon>0$ we have
$$
\mathfrak{S}_{\delta}(x; D, 2D)
\ll_{\varepsilon} (x^2 D^7)^{1/12} x^{\varepsilon}
$$
uniformly for $x\ge 3$ and $x^{6/13}\le D\le x^{2/3}$.
The same proof allows us to derive that for any $\varepsilon>0$ we have
\begin{equation}\label{UB:S(D)-1}
\mathfrak{S}_{\delta}(x; D, D')
\ll_{\varepsilon} (x^2 D^7)^{1/12} x^{\varepsilon}
\end{equation}
uniformly for $x\ge 3$, $x^{6/13}\le D\le x^{2/3}$ and $D<D'\le 2D$.
Since we have trivially
$$
\sum_{D<p^{\nu}\le D', \, \nu\ge 2} \Lambda(p^{\nu}) \psi\Big(\frac{x}{p^{\nu}+\delta}\Big)
\ll \sum_{p\le (2D)^{1/2}} \sum_{\nu\le (\log 2D)/\log p} \log p
\ll D^{1/2},
$$
the inequality \eqref{UB:S(D)-1} implies that the bound
\begin{equation}\label{UB:S(D)-2}
\sum_{D<p\le D'} (\log p) \psi\Big(\frac{x}{p+\delta}\Big)
\ll_{\varepsilon} (x^2 D^7)^{1/12} x^{\varepsilon}.
\end{equation}
holds uniformly for $x\ge 3$, $x^{6/13}\le D\le x^{2/3}$ and $D<D'\le 2D$.
Using \eqref{UB:S(D)-2}, we derive that
\begin{equation}\label{UB:S(D)-3}
\begin{aligned}
G_{2, 2}^{\langle\delta\rangle}(x)
& \ll_{\varepsilon} \max_{x^{1/2}<D\le x^{10/19}} \sum_{D<p\le 2D} \psi\Big(\frac{x}{p+\delta}\Big)
\\
& \ll_{\varepsilon} \max_{x^{1/2}<D\le x^{10/19}} 
\int_D^{2D} \frac{1}{\log t} \d \Big(\sum_{D<p\le t} (\log p) \psi\Big(\frac{x}{p+\delta}\Big)\Big)
\\
& \ll_{\varepsilon} \max_{x^{1/2}<D\le x^{10/19}} (x^2 D^7)^{1/12} x^{\varepsilon}
\\\noalign{\vskip 1,5mm}
& \ll_{\varepsilon} x^{9/19+\varepsilon}.
\end{aligned}
\end{equation}
Inserting \eqref{G21} and \eqref{UB:S(D)-3} into \eqref{G2:decomposition}, we find that
\begin{equation}\label{G2}
G_2(x) 
= \int_2^{\sqrt{x}} \frac{\d t}{\log(x/t)}
+ O\big(\sqrt{x}\exp(-c'(\log x)^{3/5}(\log_2x)^{-1/5}\big).
\end{equation}
Now the required result \eqref{thm1:eq1} follows from \eqref{G:expression}, \eqref{G1} and \eqref{G2}.

The second assertion is an immediate consequence of the first one thanks to a simple partial integration.
\hfill
$\square$

\vskip 8mm

\section{Proof of Theorem \ref{thm2}}

We begin by following the argument of \cite{LiuWuYang2021b}.
Let $f=\mathbb{1}_{\P}$ or $\mathbb{1}_{{\P}_{\rm ower}}$
and let $N\in [x^{1/3}, x^{1/2})$ be a parameter which can be chosen later.
First we write
\begin{equation}\label{4.1}
S_f(x) = \sum_{n\le x} f\Big(\Big[\frac{x}{n}\Big]\Big) = S_f^{\dagger}(x)+S_f^{\sharp}(x)
\end{equation}
with
$$
S_f^{\dagger}(x):=\sum_{n\le N} f\Big(\Big[\frac{x}{n}\Big]\Big),      
\qquad 
S_f^{\sharp}(x):=\sum_{N<n\le x} f\Big(\Big[\frac{x}{n}\Big]\Big).
$$
We have trivially
\begin{equation}\label{4.2}
S_f^{\dagger}(x)
\ll N.
\end{equation}
In order to bound $S_f^{\sharp}(x)$, we put $d=[x/n]$.
Noticing that
$$
x/n-1<d\le x/n
\;\Leftrightarrow\;
x/(d+1)<n\le x/d,
$$
we can derive that
\begin{equation}\label{4.3}
\begin{aligned}
S_f^{\sharp}(x)
& = \sum_{d\le x/N} f(d) \sum_{x/(d+1)<n\le x/d} 1
\\
& = \sum_{d\le x/N} f(d) \Big(\frac{x}{d}-\psi\Big(\frac{x}{d}\Big)-\frac{x}{d+1}+\psi\Big(\frac{x}{d+1} \Big) \Big)
\\
& = x \sum_{d\ge 1} \frac{f(d)}{d(d+1)} + \mathcal{R}_1^{f}(x, N) - \mathcal{R}_0^{f}(x, N) + O(N),
\end{aligned}
\end{equation}
where we have used the following bounds
$$
x \sum_{d>x/N} \frac{f(d)}{d(d+1)}\ll N,
\qquad
\sum_{d\le N} f(d)\Big(\psi\Big(\frac{x}{d+1}\Big) - \psi\Big(\frac{x}{d}\Big)\Big)
\ll N
$$
and
$$
\mathcal{R}_{\delta}^{f}(x, N)
= \sum_{N<d\le x/N} f(d) \psi\Big(\frac{x}{d+\delta}\Big).
$$
Combining \eqref{4.1}, \eqref{4.2} and \eqref{4.3}, it follows that
$$
S_f(x) = x \sum_{d\ge 1} \frac{f(d)}{d(d+1)} 
+ O_{\varepsilon}\big(|\mathcal{R}_1^{f}(x, N)| + |\mathcal{R}_0^{f}(x, N)| + N\big).
$$
On the other hand, we have
$$
\mathcal{R}_{\delta}^{\mathbb{1}_{{\P}_{\rm ower}}}(x, N)
= \sum_{N<p^{\nu}\le x/N} \psi\Big(\frac{x}{p^{\nu}+\delta}\Big)
= \mathcal{R}_{\delta}^{\mathbb{1}_{\P}}(x, N) + O((x/N)^{1/2}).
$$
Thus in order to prove Theorem \ref{thm2}, it suffices to show that
\begin{equation}\label{UB:Rdeltaf}
\mathcal{R}_{\delta}^{\mathbb{1}_{\P}}(x, N)\ll_{\varepsilon} N x^{\varepsilon}
\qquad
(x\ge 1)
\end{equation}
for $N=x^{9/19}$.
This can be done exactly as \eqref{UB:S(D)-3} by using \eqref{UB:S(D)-2}:
\begin{align*}
\mathcal{R}_{\delta}^{\mathbb{1}_{\P}}(x, N)
& \ll_{\varepsilon} x^{\varepsilon} \max_{x^{9/19}<D\le x^{10/19}} \sum_{D<p\le 2D} \psi\Big(\frac{x}{p+\delta}\Big)
\\
& \ll_{\varepsilon} x^{\varepsilon} \max_{x^{9/19}<D\le x^{10/19}} 
\int_D^{2D} \frac{1}{\log t} \d \Big(\sum_{D<p\le t} (\log p) \psi\Big(\frac{x}{p+\delta}\Big)\Big)
\\
& \ll_{\varepsilon} \max_{x^{9/19}<D\le x^{10/19}} (x^2 D^7)^{1/12} x^{\varepsilon}
\\\noalign{\vskip 1,5mm}
& \ll_{\varepsilon} x^{9/19+\varepsilon}.
\end{align*}
This completes the proof.
\hfill
$\square$

\vskip 3mm

\noindent{\bf Acknowledgement}.
This work is supported in part by the National Natural Science Foundation of China 
(Grant Nos. 11771211, 11971370 and 12071375).

\end{document}